\documentclass[oneside,a4paper,english,links]{amca}

\usepackage{graphicx}
\usepackage{subfigure}
\usepackage{amsmath,amsfonts}
\usepackage{commath}
\usepackage{dsfont}
\usepackage{cases}

\title{Effect of an Attached End Mass in the Dynamics of Uncertainty Nonlinear Continuous Random System}

\author[\voidaffil]{Americo Barbosa da Cunha Jr}
\author[\voidaffil]{Rubens Sampaio}

\affil[\voidaffil]{Department of Mechanical Engineering, PUC--Rio\newline
Rua Marqu\^{e}s de S\~{a}o Vicente, 225, G\'{a}vea, Rio de Janeiro - RJ, Brazil}

\everymath={\displaystyle}

\renewcommand{\vec}[1]{\ensuremath{\mathbf{#1}}}

\newcommand{\mtx}[1]{\ensuremath{\left[#1\right]}}

\newcommand{\randvar}[1]{\ensuremath{#1}}

\newcommand{\randproc}[1]{\ensuremath{#1}}

\newcommand{\pdf}[1]{\ensuremath{p_{\tiny{#1}}}}

\newcommand{\massop}[2]{\ensuremath{\int_{0}^{L} \left( \rho A \ddot{#1} (x,t) #2(x) \right) dx + m \ddot{#1} (L,t) #2(L)}}
\newcommand{\M}{\ensuremath{\mathcal{M}}}

\newcommand{\amassop}[2]{\ensuremath{ \int_{0}^{L} \rho A #1 (x,t) #2(x) dx }}
\newcommand{\aM}{\ensuremath{\widetilde{\mathcal{M}}}}

\newcommand{\dampop}[2]{\ensuremath{ \int_{0}^{L} c \dot{#1} (x,t) #2(x) dx }}
\newcommand{\C}{\ensuremath{\mathcal{C}}}

\newcommand{\stiffop}[2]{\ensuremath{\int_{0}^{L} \left( E A #1' (x,t) #2' (x) \right) dx + k #1(L,t) #2(L)}}
\newcommand{\K}{\ensuremath{\mathcal{K}}}

\newcommand{\forceop}[1]{\ensuremath{ \int_{0}^{L} f(x,t) #1 (x) dx }}
\newcommand{\F}{\ensuremath{\mathcal{F}}}

\newcommand{\nlforceop}[2]{\ensuremath{ -k_{NL} \left[ u(L,t) \right]^3 #2 (L) }}
\newcommand{\NL}{\ensuremath{\mathcal{F}_{NL}}}

\newcommand{\indfunc}[1]{\ensuremath{\mathds{1}_{#1} }}

\newcommand{\entropy}[1]{\ensuremath{\mathbb{S}\left[  #1 \right] }}

\newcommand{\entropyop}[3]{\ensuremath{ - \int_{#2}^{#3} #1(\xi) \ln \left[ #1(\xi) \right] d\xi}}

\newcommand{\expectation}[1]{\ensuremath{\mathbb{E}\left[  #1 \right] }}

\newcommand{\expectationop}[3]{\ensuremath{\int_{#2}^{#3} #1(\xi) \pdf{#1}(\xi) d\xi} }

\newcommand{\R}{\ensuremath{\mathbb{R}}}           
\newcommand{\U}{\ensuremath{\mathcal{U} }}   
\newcommand{\W}{\ensuremath{\mathcal{W} }}   

\renewcommand{\P}[1][]{\ensuremath{{\mathbb{P}^{#1}} }}
\newcommand{\A}[1][]{\ensuremath{{\mathbb{A}^{#1}} }}

\begin{document}
\vspace{3cm}

\maketitle


\begin{abstract}

This work studies the dynamics of a one dimensional elastic bar with 
random elastic modulus and prescribed boundary conditions, say, fixed at one end, 
and attached to a lumped mass and two springs (one linear and another nonlinear) 
on the other extreme. The system analysis assumes that the elastic modulus has 
gamma probability distribution and uses Monte Carlo simulations to compute 
the propagation of uncertainty in this continuous--discrete system.
After describing the deterministic and the stochastic modeling of the system, 
some configurations of the model are analyzed in order to characterize the effect 
of the lumped mass in the overall behavior of this dynamical system.

\end{abstract}

\section{INTRODUCTION}

The dynamics of a mechanical system depends on some parameters such as
physical and geometrical properties, constraints, external and internal loading, 
initial and boundary conditions.
Most of the theoretical models used to describe the behavior of
a mechanical system assume nominal values for these parameters,
such that the model gives one response for a given particular input.
In this case the system is \emph{deterministic} and its behavior is 
described by a single set of differential equations.
However, in real systems they do not have a fixed value
since they are subjected to uncertainties of measurement, 
imperfections in manufacturing processes, change of properties, etc. 
This variability in the set of system parameters leads to a
large number of possible system responses for a given particular input. 
Now the system is \emph{stochastic} and there is a family of 
differential equations sets (one for each realization of the random system)
associated to it.

This work aims to study the propagation of uncertainty
in the dynamics of a nonlinear continuos random system with a discrete
element attached to it. In this sense, this work considers a one dimensional 
elastic bar, with random elastic modulus, fixed on the left extreme and 
with a lumped mass and two springs (one linear and another nonlinear) 
on the right extreme (fixed-mass-spring bar).

This paper is organized as follows. In section~\ref{determ_approach}
is presented the deterministic modeling of the problem, the discretization 
procedure and the algorithm used to solve the equation of interest.
The stochastic modeling of the problem is shown in 
section~\ref{stochastic_approach}, as well as the construction of a 
probability distribution for the elastic modulus, using the 
maximum entropy principle, and a brief discussion on the Monte Carlo method. 
In section~\ref{num_experim}, some configurations of the model are analyzed 
in order to characterize the effect of lumped mass in the system dynamical behavior. 
Finally, in section~\ref{concl_remaks}, the main conclusions are emphasized 
and some directions for future work outlined.


\section{DETERMINISTIC APPROACH}
\label{determ_approach}

\begin{figure}[h!]
	\centering
	\includegraphics[scale=0.6]{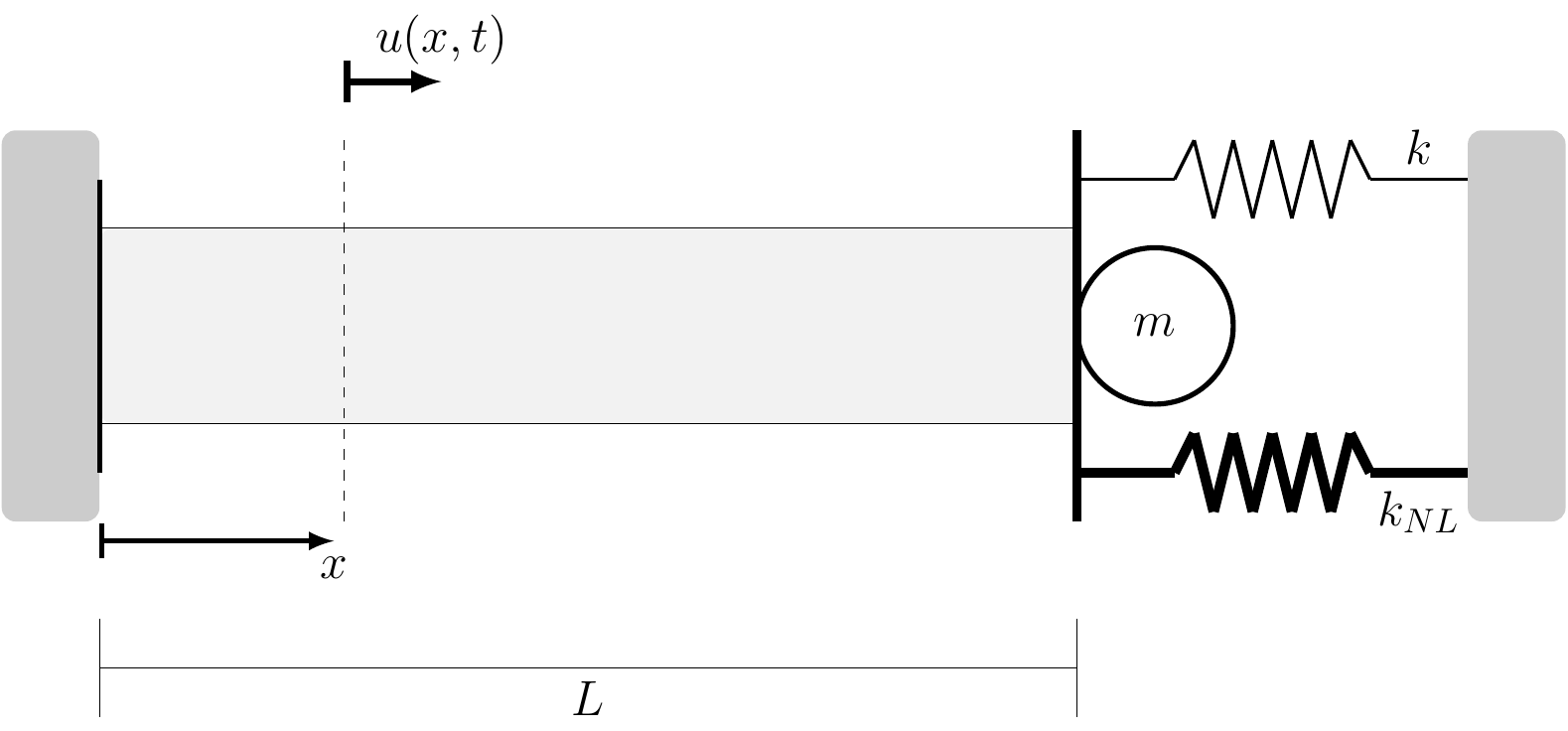}
	\caption{Sketch of a bar fixed at one and attached
					to two springs and a lumped mass on the other extreme.}
	\label{bar_fig}
\end{figure}

The continuous system of interest is the one-dimensional 
fixed-mass-spring bar shown in Figure~\ref{bar_fig}. 

\subsection{Strong Formulation}
\label{strong_form}

The displacement of this system $u$ evolves according 
to the following partial differential equation

\begin{equation}
    \rho A \dpd[2]{u}{t}(x,t) + c \dpd{u}{t}(x,t) = 
    \dpd{}{x}\left( E A \dpd{u}{x}(x,t)  \right)  + f(x,t),
    \label{bar_eq}
\end{equation}

\noindent
which is valid for $0 < x < L$ and $0 < t < T$, being 
$L$ the bar unstretched length and $T$ a finite instant of time.
In this equation $\rho$ is the mass density, $E$ is the elastic modulus, 
$A$ is the circular cross section area, $c$ is the damping coefficient, and
$f(x,t)$ is an external force depending on position $x$ and instant $t$.

The left side of the bar is fixed at a rigid wall while the right 
side is attached to a lumped mass $m$ and two springs 
fixed to a rigid wall. The first spring (of stiffness $k$) 
is linear and exerts a restoring force proportional 
to the stretching on the bar. The second spring 
(of stiffness $k_{NL}$) is nonlinear and its restoring force is 
proportional to the cube of the stretching. The force which 
the lumped mass exerts on the bar is proportional to acceleration.
These boundary conditions read as

\begin{equation}
    u(0,t) = 0 \qquad \mbox{and} \qquad 
    E A \dpd{u}{x} (L,t) = -k u(L,t) - k_{NL} \left[ u(L,t)\right]^3 - m \dpd[2]{u}{t}(L,t).
    \label{bar_bc}
\end{equation}

Initially, any point $x$ of the bar presents displacement and a velocity respectively
equal to

\begin{equation}
    u(x,0) = u_0(x) \qquad \mbox{and} \qquad \dpd{u}{t}(x,0) = \dot{u}_0(x),
    \label{bar_ic}
\end{equation}

\noindent
for $0 \leq x \leq L$. In these equations $u_0$ and $\dot{u}_0$ are given
functions of position $x$.

Moreover, it is noteworthy that $u$ is assumed to be as regular as needed 
for the initial--boundary value problem of Eqs.(\ref{bar_eq}), (\ref{bar_bc}), 
and (\ref{bar_ic}) to be well posed.

\subsection{Variational Formulation}
\label{var_form}

Let $\U_t$ be the class of (time dependent) basis functions
and $\W$ be the class of weight functions. These sets are 
chosen as the space of functions with square integrable spatial derivative,
which satisfy the essential boundary condition defined by Eq.(\ref{bar_bc}).

The variational formulation of the problem under study says that one wants 
to find $u \in \U_t$ that satisfy, for all $w \in \W$, the weak equation of
motion given by

\begin{equation}
    \M(\ddot{u},w) + \C(\dot{u},w) + \K(u,w) = \F(w) + \NL(u,w),
    \label{weak_eq}
\end{equation}

\noindent
where 
$\M$ is the mass operator,
$\C$ is the damping operator,
$\K$ is the stiffness operator, 
$\F$ is the external force operator, and
$\NL$ is the nonlinear force operator.
These operators are, respectively, defined as


\begin{equation}
    \M(\ddot{u},w) = \massop{u}{w}, \label{mass_op}
\end{equation}    

\begin{equation}    
    \C(\dot{u},w) = \dampop{u}{w}, \label{damp_op}
\end{equation}

\begin{equation}
    \K(u,w) = \stiffop{u}{w}, \label{stiff_op}
\end{equation}

\begin{equation}
    \F(w) = \forceop{w}, \label{force_op}
\end{equation}

\begin{equation}
    \NL(u,w) = \nlforceop{u}{w}, \label{nlforce_op}
\end{equation}

\noindent
where ~$\dot{}$~ is an abbreviation for temporal derivative and
~$'$~ is an abbreviation for spatial derivative.

The variational formulations for the initial conditions of Eq.(\ref{bar_ic}),
which are valid for all $w \in \W$, are respectively given by

\begin{equation}
	\aM(u(\cdot,0),w) = \aM(u_0,w),
    \label{weak_ic_eq1}
\end{equation}

\noindent
and

\begin{equation}
	\aM(\dot{u}(\cdot,0),w) = \aM(\dot{u}_0,w),
    \label{weak_ic_eq2}
\end{equation}

\noindent
where $\aM$ is the associated mass operator, defined as

\begin{equation}
    \aM(u,w) = \amassop{u}{w}.
    \label{amass_op}
\end{equation}

\subsection{An Eigenvalue Problem}
\label{eigen_probl}

Now consider the following generalized eigenvalue problem
associated to Eq.(\ref{weak_eq}),

\begin{equation}
		- \nu^2 \M(\phi,w) + \K(\phi,w) = 0,
		\label{gen_eig_eq}
\end{equation}

\noindent
where $\nu$ is a natural frequency and $\phi$ is an associated mode shape.

In order to solve Eq.(\ref{gen_eig_eq}), the technique of 
separation of variables is employed, which leads to a 
Sturm-Liouville problem \citep{algwaiz2007},
with denumerable number of solutions.
Therefore, this generalized eigenvalue problem
has a denumerable number of solutions,
all of then such as the following eigenpair
$(\nu_n^2,\phi_n)$, where $\nu_n$ is the
$n$-th bar natural frequency and
$\phi_n$ is the $n$-th bar mode shape.

It is important to observe that, the eigenfunctions
$\{\phi_n\}_{n=1}^{+\infty}$ span the space of 
functions which contains the solution of the 
Eq.(\ref{gen_eig_eq})  \citep{brezis2010}.
As can be seen in \cite{hagedorn2007}, 
these eigenfunctions satisfy, for all  $m \neq n$,
the orthogonality relations given by

\begin{equation}
		\M(\phi_n,\phi_m) = 0,
		\label{modalM_eq}
\end{equation}

\noindent
and

\begin{equation}
		\K(\phi_n,\phi_m) = 0,
		\label{modalK_eq}
\end{equation}

\noindent
which made then good choices for the basis function when a weighted 
residual procedure \citep{scriven1966p735} is used to approximate 
the solution of a nonlinear variational equation, 
such as Eq.(\ref{weak_eq}).

\subsection{Mode Shapes and Natural Frequencies}
\label{mode_shapes_nat_freq}

According to \cite{blevins1993}, a fixed-mass-spring bar 
has its natural frequencies and the corresponding orthogonal modes shape given by

\begin{equation}
		\nu_n = \lambda_n \frac{\bar{c}}{L},
		\label{nat_freq_eq}
\end{equation}

\noindent
and

\begin{equation}
		\phi_n (x) = \sin{ \left( \lambda_n \frac{x}{L} \right) },
		\label{shape_mod_eq}
\end{equation}

\noindent
where $\bar{c} = \sqrt{E/\rho}$ is the wave speed, and the 
$\lambda_n$ are the solutions of

\begin{equation}
		\cot{ \left( \lambda_n \right) }  + 
		\left( \frac{kL}{AE} \right) \frac{1}{\lambda_n} -
		\left( \frac{m}{\rho AL} \right) \lambda_n = 0.
		\label{lambda_eq}
\end{equation}

\begin{figure} [h!]
				\centering
				\subfigure[$\nu_1 \approx 0.07 \times 10^5~rad/s$]{
				\includegraphics[scale=0.28]{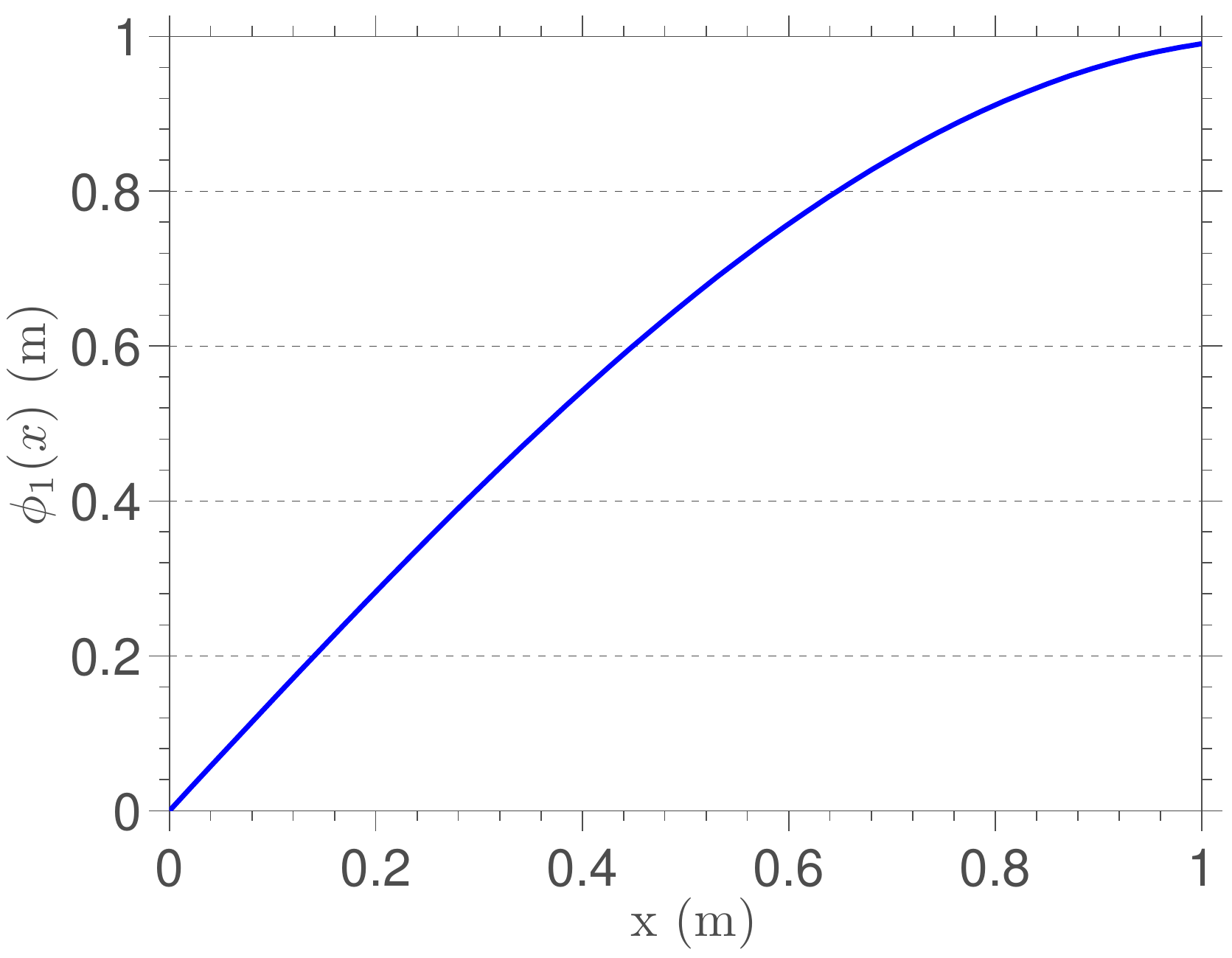}}
				\subfigure[$\nu_2 \approx 0.22 \times 10^5~rad/s$]{
				\includegraphics[scale=0.28]{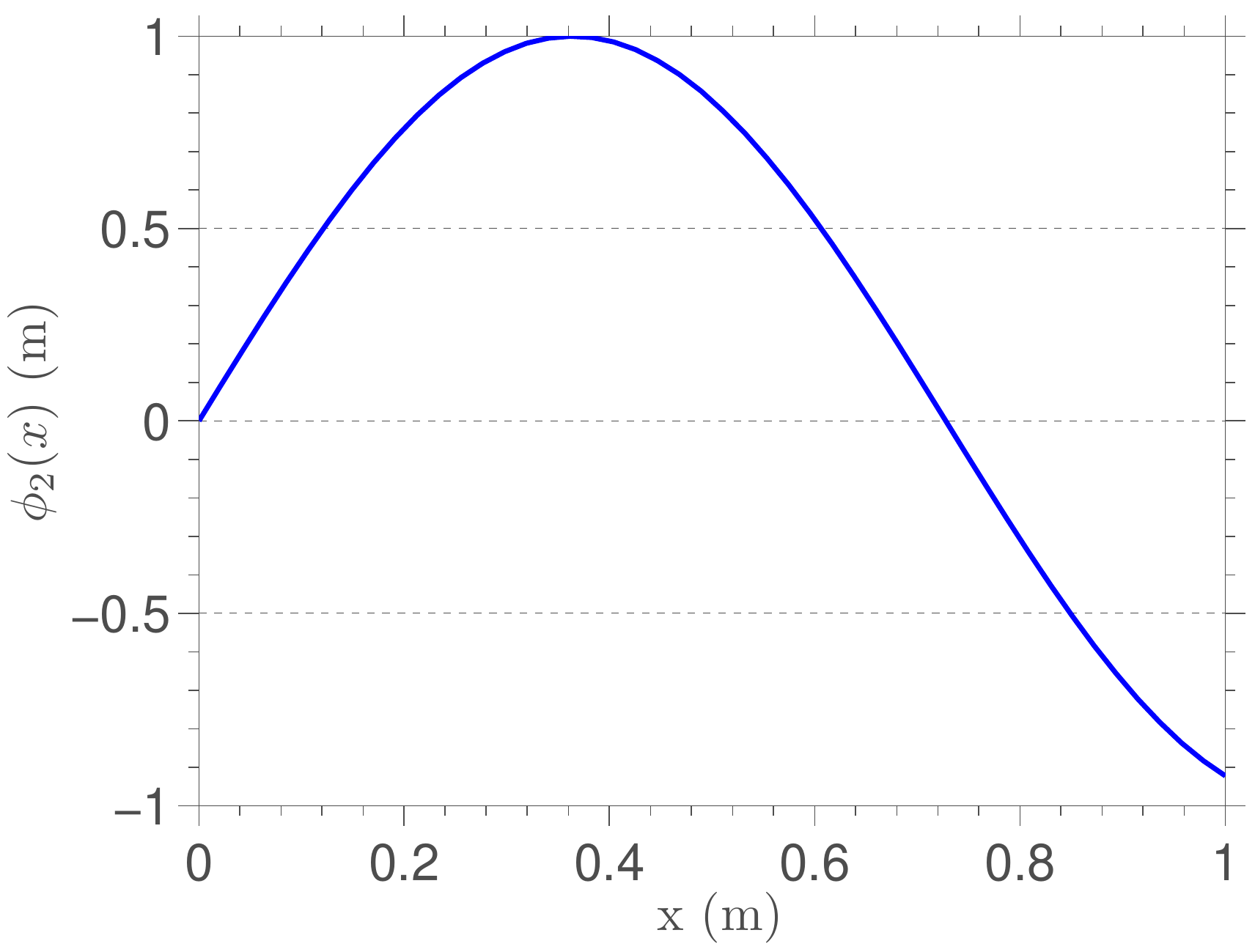}}
				\subfigure[$\nu_3 \approx 0.37 \times 10^5~rad/s$]{
				\includegraphics[scale=0.28]{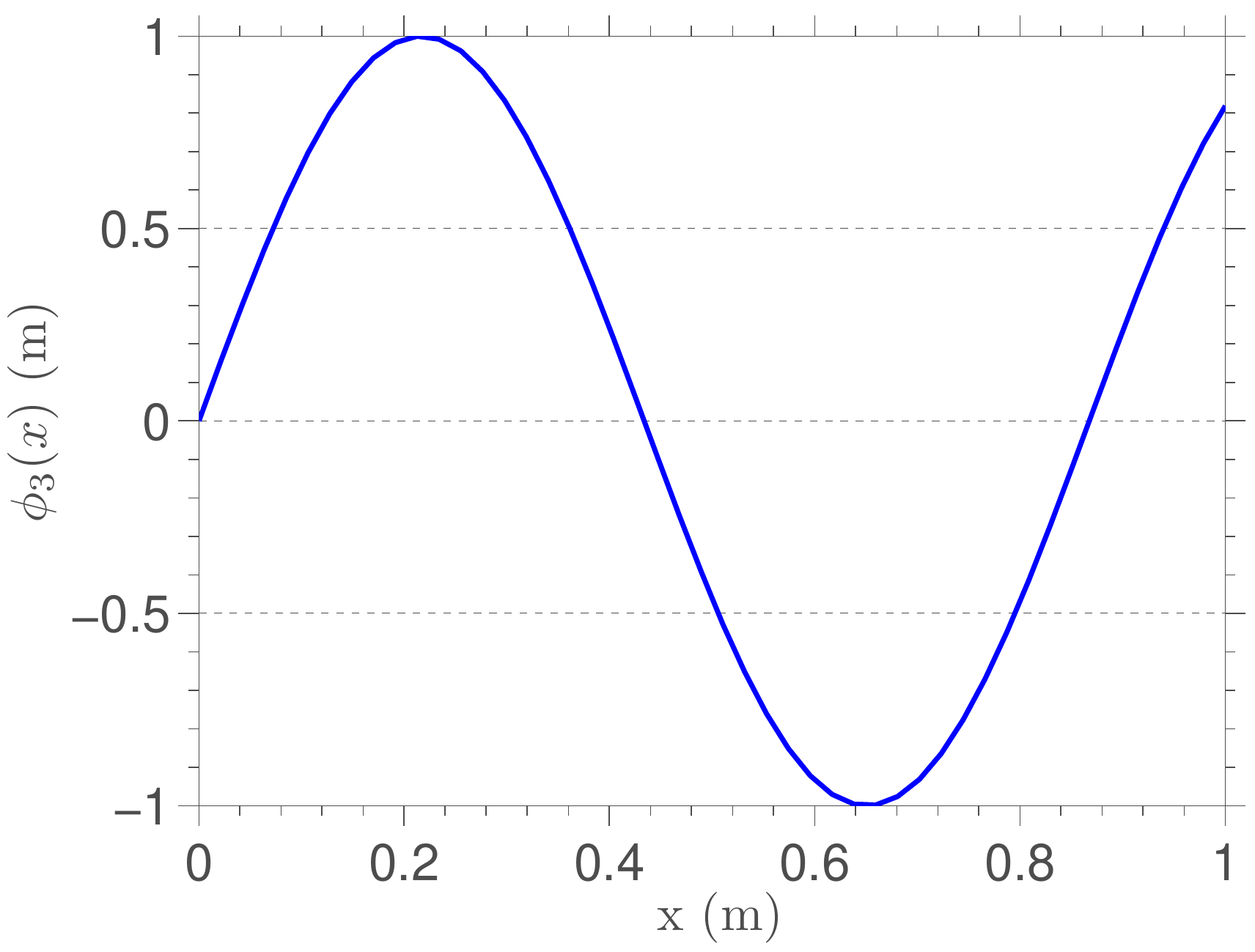}}\\
				\subfigure[$\nu_4 \approx 0.52 \times 10^5~rad/s$]{
				\includegraphics[scale=0.28]{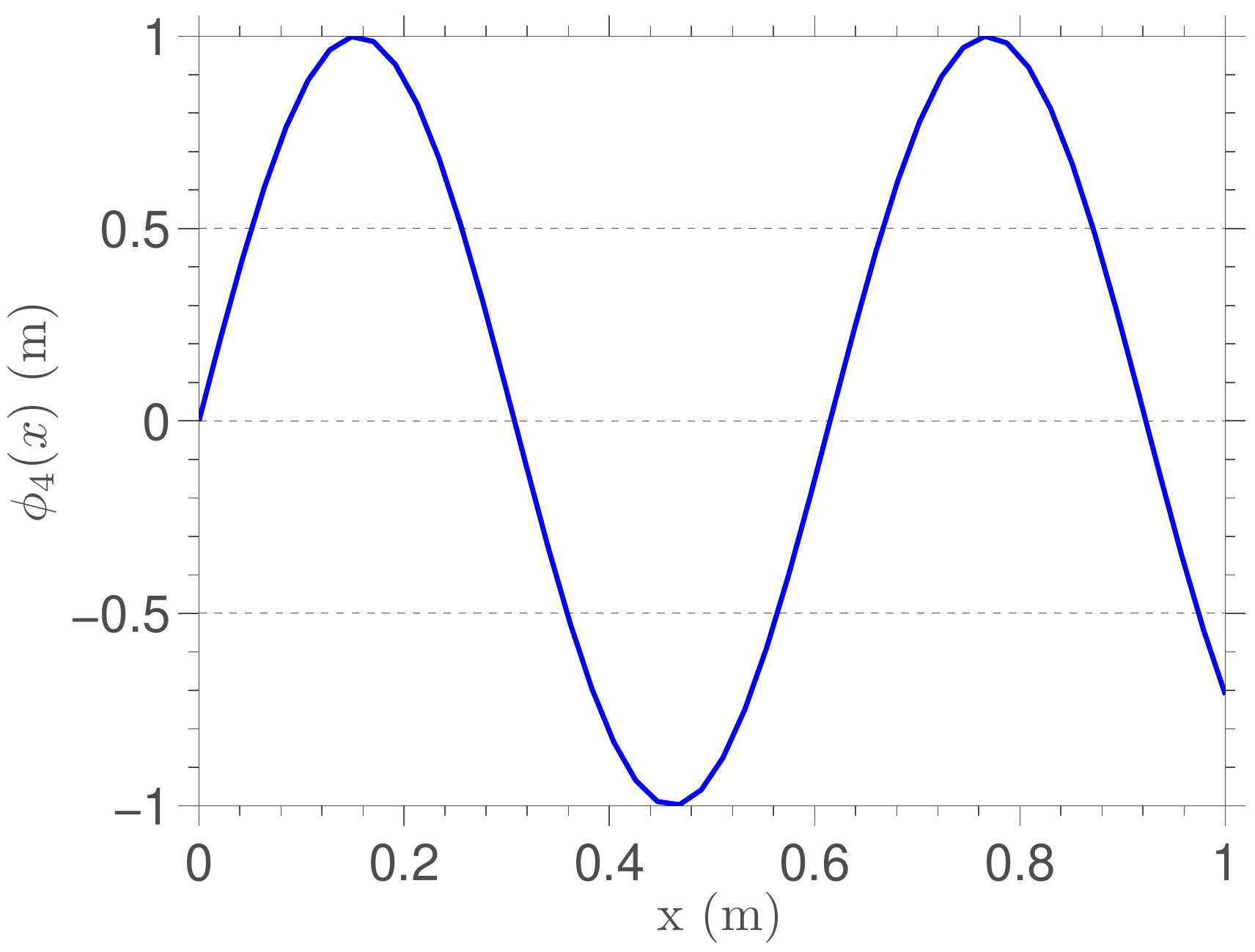}}
				\subfigure[$\nu_5 \approx 0.67 \times 10^5~rad/s$]{
				\includegraphics[scale=0.28]{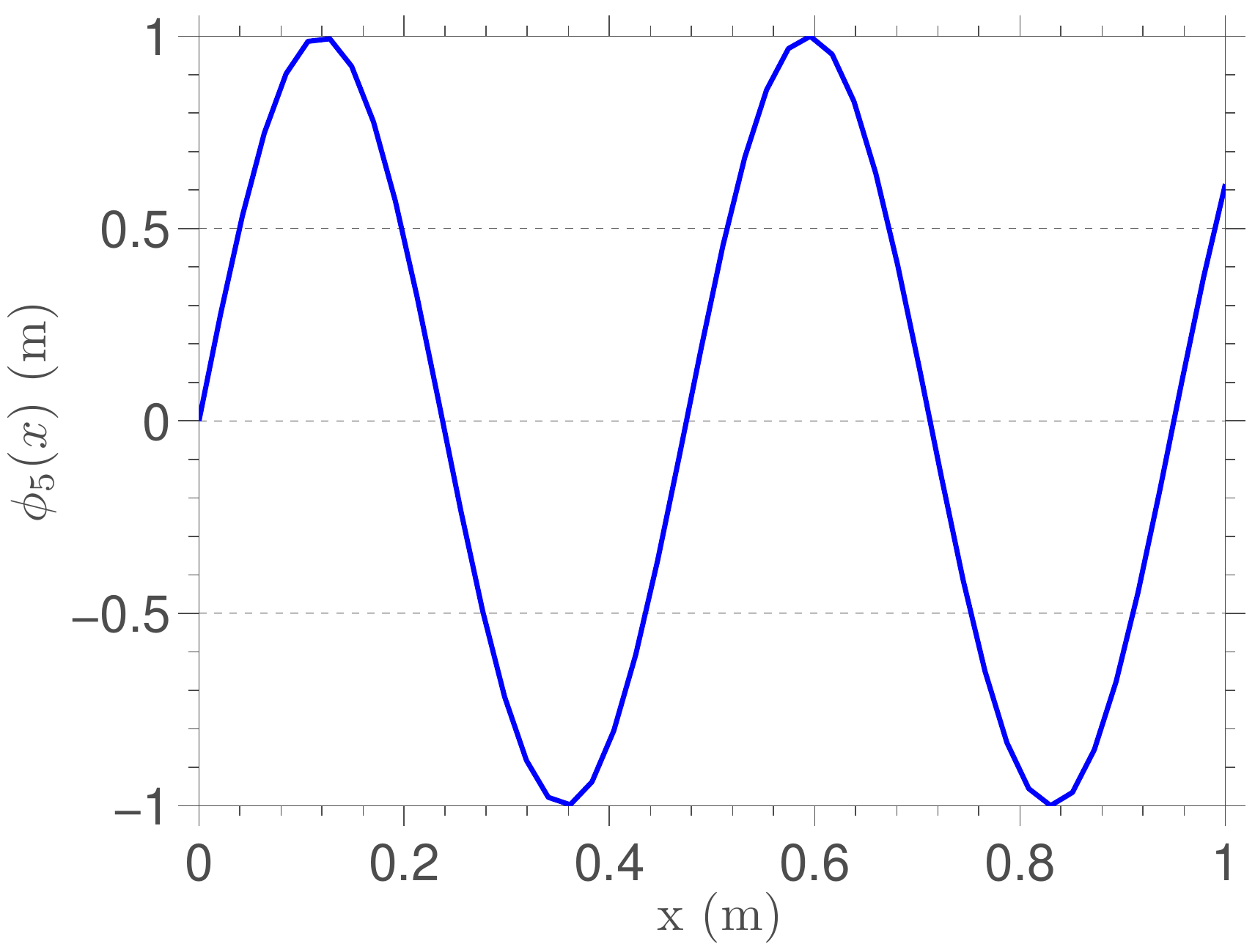}}
				\subfigure[$\nu_6 \approx 0.83 \times 10^5~rad/s$]{
				\includegraphics[scale=0.28]{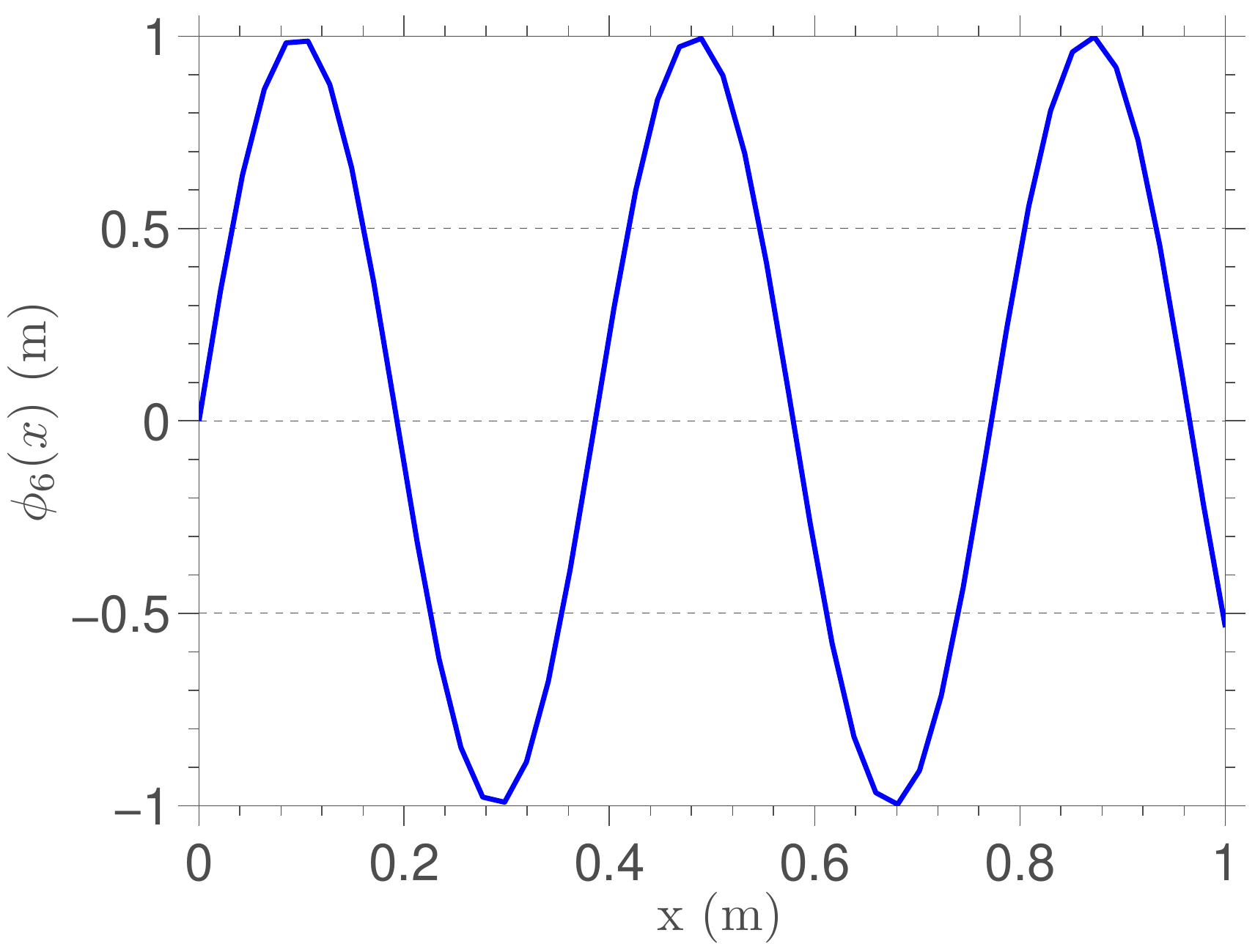}}
				\caption{The first six orthogonal modes shape and the corresponding 
				(approximated) natural frequencies of a fixed-mass-spring bar with $m=1.5~kg$.}
				\label{shape_modes_fig}
\end{figure}

The first six orthogonal modes shape of the fixed-mass-spring bar with $m=1.5~kg$,
whose the other parameters are presented in the beginning of section \ref{num_experim},
are illustrated in Figure~\ref{shape_modes_fig}. In this figure each sub-caption 
indicates the approximated natural frequency associated with the 
corresponding mode.

\subsection{Galerkin Formulation}
\label{galerkin_form}

In order of approximate the solution of 
Eqs.(\ref{weak_eq}),~(\ref{weak_ic_eq1})~and~(\ref{weak_ic_eq2})
the Galerkin method \citep{hughes2000} is employed. 
Therefore, the displacement field $u$ is approximated by a linear 
combination of the form

\begin{equation}
    u^N(x,t) = \sum_{n=1}^N u_n(t) \phi_n(x),
    \label{uN_eq}
\end{equation}

\noindent
where the basis functions $\phi_n$ are the orthogonal modes shape 
of the fixed-mass-spring bar, exemplified in the end of
section~\ref{mode_shapes_nat_freq}, and the coefficients $u_n$ 
are time-dependent functions. For a reason that
will be clear soon, define $\vec{u}(t)$ of $\R^N$
as the vector in which the $n$-th component is $u_n(t)$.

Since $u^N$ is not a solution of Eq.(\ref{weak_eq}), when the field 
$u$ is approximated by $u^N$ a residual function is obtained. 
This residual function is orthogonally projected into the vector 
space spanned by the functions $\{\phi_n\}_{n = 1}^{N}$
in order to minimize the error incurred by the approximation 
\citep{hughes2000}. This procedure results in the following 
$N \times N$ set of nonlinear ordinary differential equations 

\begin{equation}
    \mtx{M} \ddot{\vec{u}}(t) + \mtx{C} \dot{\vec{u}}(t) +
    \mtx{K}\vec{u}(t) = \vec{f}(t) + \vec{f}_{NL}\left( \vec{u}(t) \right),
    \label{galerkin_eq}
\end{equation}

\noindent
supplemented by the following pair of initial conditions

\begin{equation}
    \vec{u}(0) =  \vec{u}_0 \qquad \mbox{and} \qquad \vec{\dot{u}}(0) =  \vec{\dot{u}}_0.
    \label{galerkin_ic_eq}
\end{equation}

\noindent
where $\mtx{M}$ is the mass matrix, $\mtx{C}$ is the damping matrix, 
$\mtx{K}$ is the stiffness matrix, and the upper dot again denotes 
the time derivative. Also, $\vec{f}(t)$, $\vec{f}_{NL}\left(\vec{u}(t)\right)$, 
$\vec{u}_0$, and $\vec{\dot{u}}_0$ are vectors of $\R^N$, which respectively 
represent the external force, the nonlinear force, the initial position, 
and the initial velocity. 

The initial value problem of
Eqs.(\ref{galerkin_eq}) and (\ref{galerkin_ic_eq}) has its solution approximated 
by Newmark method \citep{newmark1959p67}. The reader interested in details
about this integration scheme is encouraged to see \cite{hughes2000}.


\section{STOCHASTIC APPROACH}
\label{stochastic_approach}

\subsection{Probabilistic Model}

Consider a probability space  $(\Omega, \A, \P)$, 
where $\Omega$ is sample space, 
$\A$ is a $\sigma$-field over $\Omega$ and 
$\P$ is a probability measure. 
In this probabilistic space, the elastic modulus 
is assumed to be a random variable $\randvar{E}: \Omega \to \R$ that 
associates to each event $\omega \in \Omega$ 
a real number $\randvar{E}(\omega)$.
Consequently, the displacement of the bar is the random field
\mbox{$\randproc{U}: [0,L] \times [0,T]  \times \Omega \to \R$,} 
which evolves according the following stochastic 
partial differential equation

\begin{equation}
    \rho A \dpd[2]{\randproc{U}}{t} (x,t,\omega) + 
    c \dpd{\randproc{U}}{t} (x,t,\omega) = 
    \dpd{}{x} \left( \randvar{E}(\omega) A \dpd{\randproc{U}}{x} (x,t,\omega)\right)
    + f(x,t),
    \label{randbar_eq}
\end{equation}

\noindent
being the partial derivatives now defined in the mean square
sense \citep{papoulis2002}. This problem has boundary and 
initial conditions similar to those defined in 
Eqs.(\ref{bar_bc})~and~(\ref{bar_ic}), by changing 
$u$ for $\randproc{U}$ only.

\subsection{Elastic Modulus Distribution}

The elastic modulus cannot be negative, so it is reasonable to assume 
the support of random variable $\randvar{E}$ as the interval $(0,+\infty)$. 
Therefore, the probability density function (PDF) of $\randvar{E}$ 
is a nonnegative function $\pdf{\randvar{E}}: (0,+\infty) \to \R$, which 
respects the following normalization condition

\begin{equation}
	\int_{0}^{+\infty} \pdf{\randvar{E}}(\xi) d\xi = 1.
	\label{norm_cond_E}
\end{equation}

Also, the mean value of $\randvar{E}$ is known real number 
$\mu_{\randvar{E}}$, i.e.,

\begin{equation}
	\expectation{\randvar{E}} = \mu_{\randvar{E}},
	\label{meanval_E}
\end{equation}

\noindent
where the expected value operator of $\randvar{E}$ is defined as

\begin{equation}
	\expectation{\randvar{E}} = \expectationop{\randvar{E}}{0}{+\infty}
	\label{def_expval_op}
\end{equation}

Finally, one also wants that $\randvar{E}$ has finite variance, i.e., 

\begin{equation}
	\expectation{\left(\randvar{E} - \mu_{\randvar{E}} \right)^2} < + \infty,
	\label{finite_var_E}
\end{equation}

\noindent
which is possible \citep{soize2000p277}, for example, if

\begin{equation}
	\expectation{\ln \left( \randvar{E} \right) } < + \infty.
	\label{finite_log_E}
\end{equation}

Following the suggestion of \cite{soize2000p277},
the maximum entropy principle
\citep{shannon1948p379,jaynes1957p620,jaynes1957p171}
is employed in order to consistently specify $\pdf{\randvar{E}}$.
This methodology chooses for $\randvar{E}$ the PDF
which maximizes the differential entropy function,
defined by

\begin{equation}
	\entropy{\pdf{\randvar{E}}} = \entropyop{\pdf{\randvar{E}}}{0}{+\infty},
	\label{dif_entropy}
\end{equation}

\noindent
subjected to (\ref{norm_cond_E}), (\ref{meanval_E}), and (\ref{finite_log_E}),
the restrictions that effectively define the known information about $\randvar{E}$.

Respecting the constraints imposed by (\ref{norm_cond_E}),
(\ref{meanval_E}), and (\ref{finite_log_E}), the PDF that maximizes 
Eq.(\ref{dif_entropy}) is given by

\begin{equation}
	\pdf{\randvar{E}}(\xi) = \indfunc{(0,+\infty)}
	\frac{1}{\mu_{\randvar{E}}} 
	\left( \frac{1}{\delta_{\randvar{E}}^2} \right)^{ \left( \displaystyle \frac{1}{\delta_{\randvar{E}}^2} \right) }
	\frac{1}{\Gamma(1/\delta_{\randvar{E}}^2)} 
	\left( \frac{\xi}{\mu_{\randvar{E}}} \right)^{ \left( \displaystyle \frac{1}{\delta_{\randvar{E}}^2}-1 \right) }
	\exp\left( - \frac{\xi}{\delta_{\randvar{E}}^2 \mu_{\randvar{E}}} \right),
	\label{gamma_distrib}
\end{equation}

\noindent
where $\indfunc{(0,+\infty)}$ denotes the indicator function of the 
interval $(0,+\infty)$, $\delta_{\randvar{E}}$ is the dispersion factor
of $\randvar{E}$, and $\Gamma$ indicates the gamma function.
This PDF is a gamma distribution. 

\subsection{Stochastic Solver: Monte Carlo Method}
\label{stochastic_solver}

Uncertainty propagation in the stochastic dynamics of the continuous--discrete 
system under study is computed by Monte Carlo (MC) method 
\citep{metropolis1949p335}. 
This stochastic solver uses a Mersenne twister pseudorandom 
number generator \citep{matsumoto1998p3}, 
to obtain many realizations of the random 
variable $\randvar{E}$. Each one of these realizations
defines a new Eq.(\ref{weak_eq}), so that a
new variational problem is obtained.
After that, these new variational problems
are solved deterministically, such as in section~\ref{galerkin_form}.
All the MC simulations reported in this work use $4^5$ samples
to access the random system. Further details about MC method 
can be seen in \cite{liu2001,shonkwiler2009,casella2010}.


\section{NUMERICAL EXPERIMENTS}
\label{num_experim}

The numerical experiments presented in this section adopt
the following deterministic parameters for the studied system:
$\rho = 7900~kg/m^3$, $c=10~kNs/m$, $A = 625\pi~mm^2$, 
$k = 650~N/m$, $k_{NL} = 650 \times 10^{13}~N/m^3$, $L= 1~m$, 
and $T=8~ms$. Besides that, four values for the lumped mass are considered:
$m=1.5,~7.5,~15,~\mbox{and}~75~kg$. The random variable $\randvar{E}$, 
is characterized by $\mu_{\randvar{E}} = 203~GPa$ and $\delta_{\randvar{E}} = 10\%$.

The initial conditions for displacement and velocity are respectively given by

\begin{equation}
    u_0 = \alpha_1 \phi_3(x) + \alpha_2 x, \qquad \mbox{and} \qquad \dot{u}_0 = 0,
\end{equation}

\noindent
where $\alpha_1 =0.1~mm$ and $\alpha_2 =0.5 \times 10^{-3}$. Note that $u_0$ reaches 
the maximum value at $x = L$. This function is used to ``activate" the spring cubic 
nonlinearity, which depends on the displacement at $x = L$.

The time-dependent external force acting on the system
has the form of a sine wave with circular frequency equal 
to the first natural frequency

\begin{equation}
	f(x,t) = \sigma \phi_1(x) \sin{ \left( \nu_{1} t \right) },
	\label{ext_force}
\end{equation}

\noindent
where the external force amplitude is $\sigma = 1~N$.

\begin{figure} [h!]
				\centering
				\subfigure[$\rho A L/m \approx 10$]{
				\includegraphics[scale=0.38]{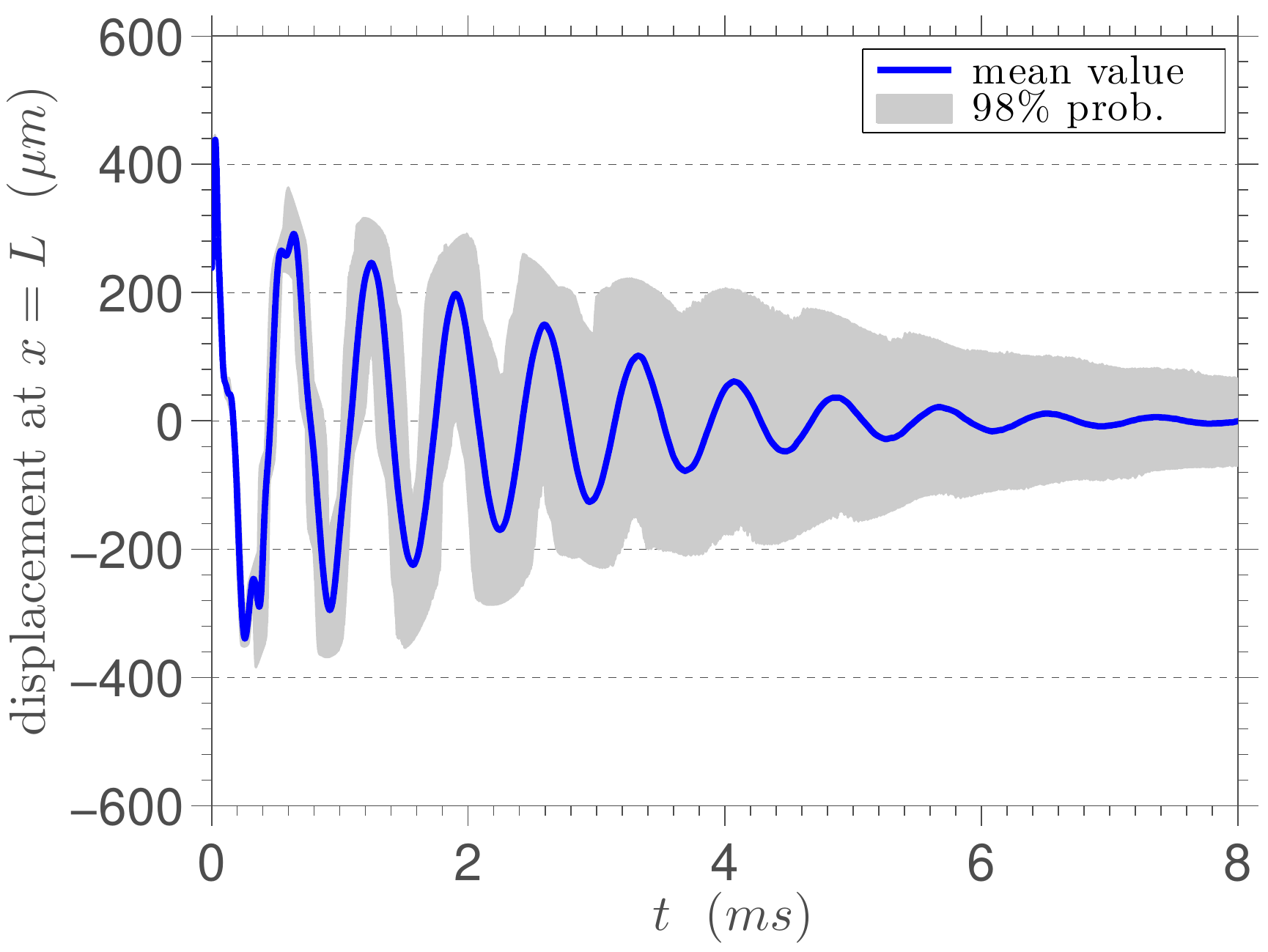}}
				\subfigure[$\rho A L/m \approx 2$]{
				\includegraphics[scale=0.38]{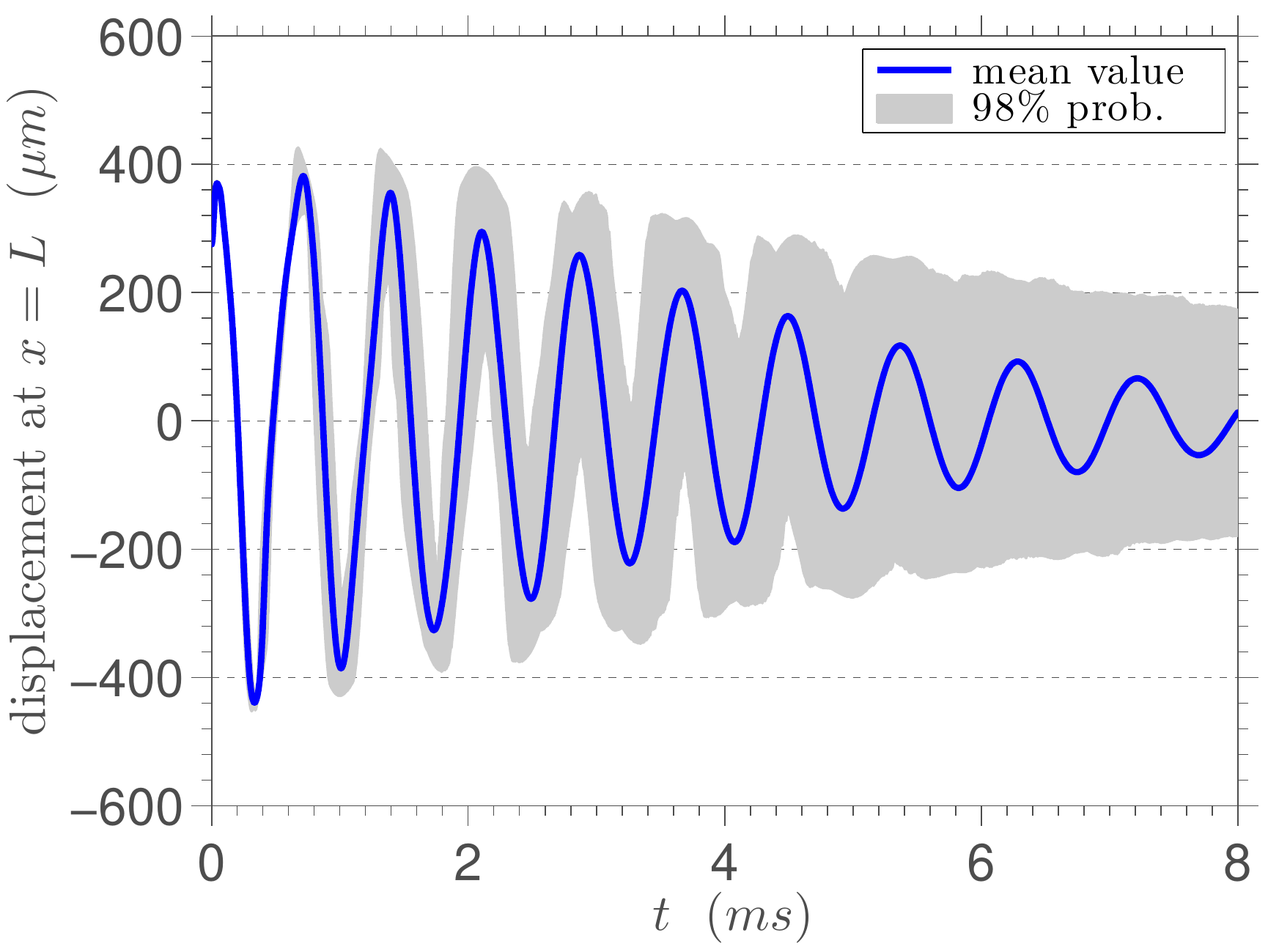}}\\
				\subfigure[$\rho A L/m \approx 1$]{
				\includegraphics[scale=0.38]{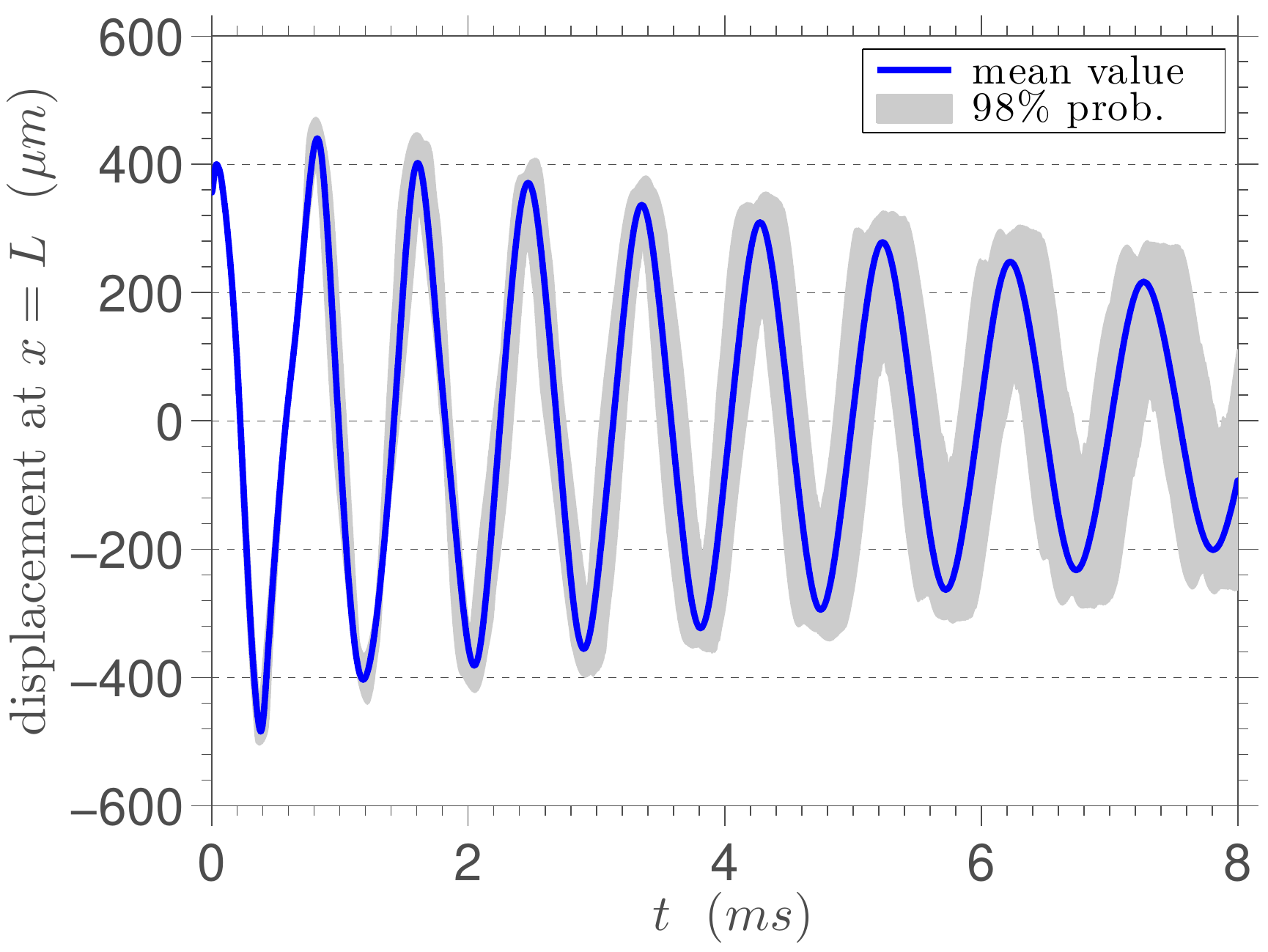}}
				\subfigure[$\rho A L/m \approx 0.2$]{
				\includegraphics[scale=0.38]{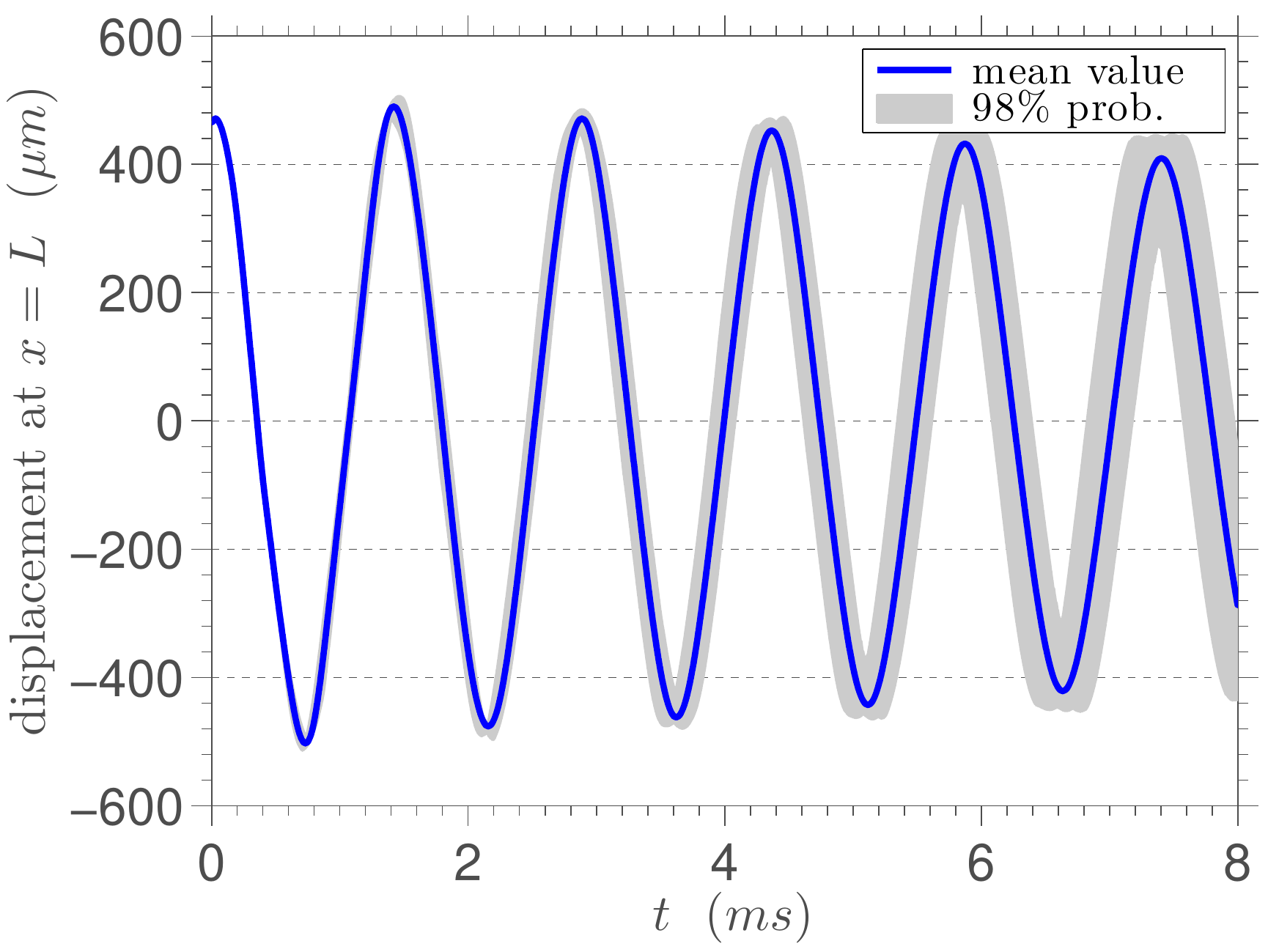}}
				\caption{This figure illustrates the mean value (blue line) and a 98\% of probability interval 
							  of confidence (grey~shadow) for the random process $\randproc{U}(L,\cdot,\cdot)$,
							  for several values of the continuous--discrete mass ratio.}
				\label{ci_uL_fig}
\end{figure}

\subsection{Analysis of Random System Envelope of Reliability and Phase Space}
\label{analysis_phase_space}

The mean value of $\randproc{U}(L,\cdot,\cdot)$ and an envelope of reliability, 
wherein a realization of the stochastic system has 98\% of probability of being contained,
are shown, for different values of the continuous--discrete mass ratio $\rho A L/m$,
in Figure~\ref{ci_uL_fig}. By observing this figure one can note that, as the value of lumped mass increases, 
the decay of the system displacement amplitude decreases significantly. This indicates 
that this continuous--discrete system is not much influenced by damping for small values of
the continuous--discrete mass ratio.

\begin{figure} [ht!]
				\centering
				\subfigure[$\rho A L/m \approx 10$]{
				\includegraphics[scale=0.38]{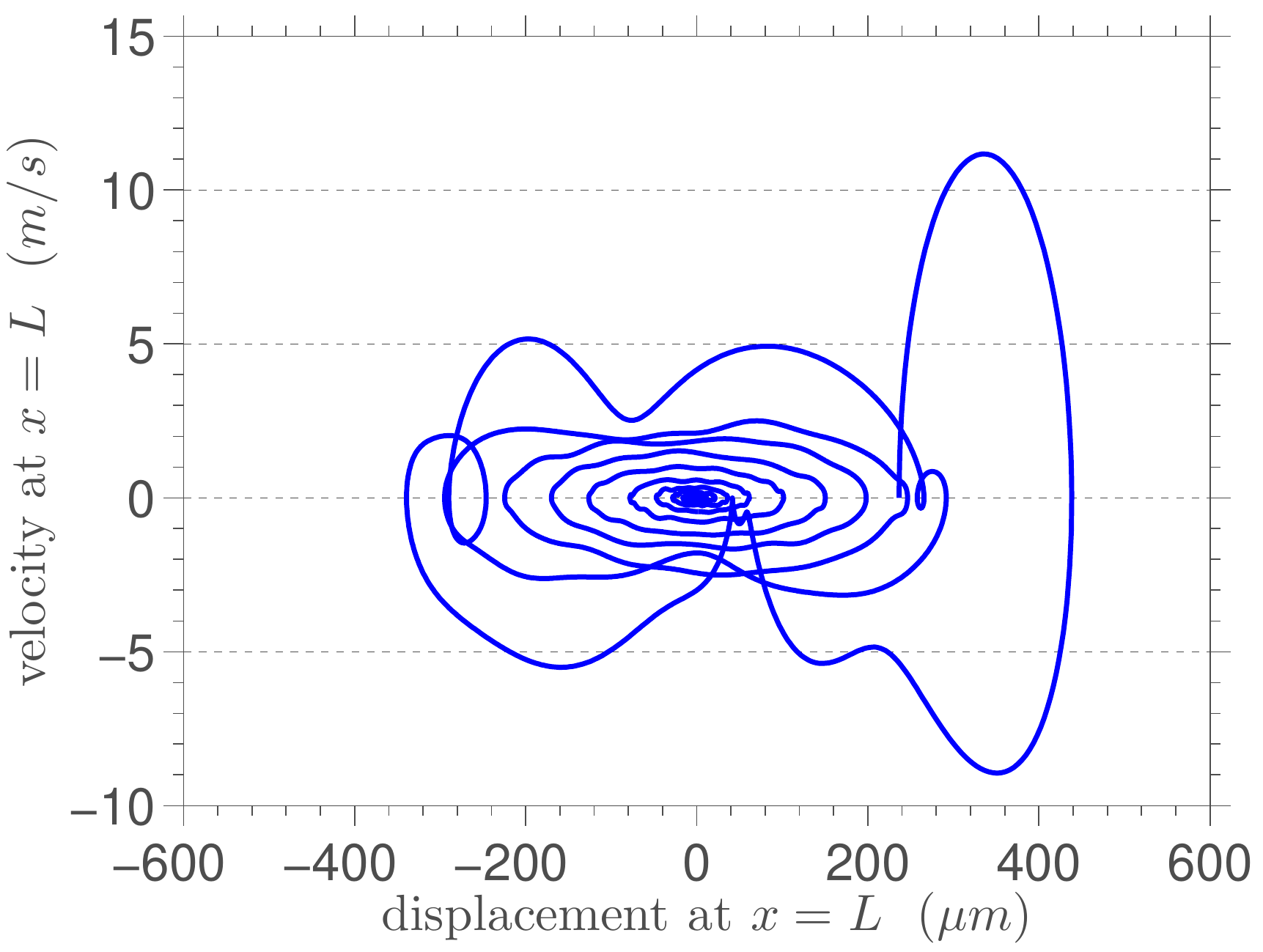}}
				\subfigure[$\rho A L/m \approx 2$]{
				\includegraphics[scale=0.38]{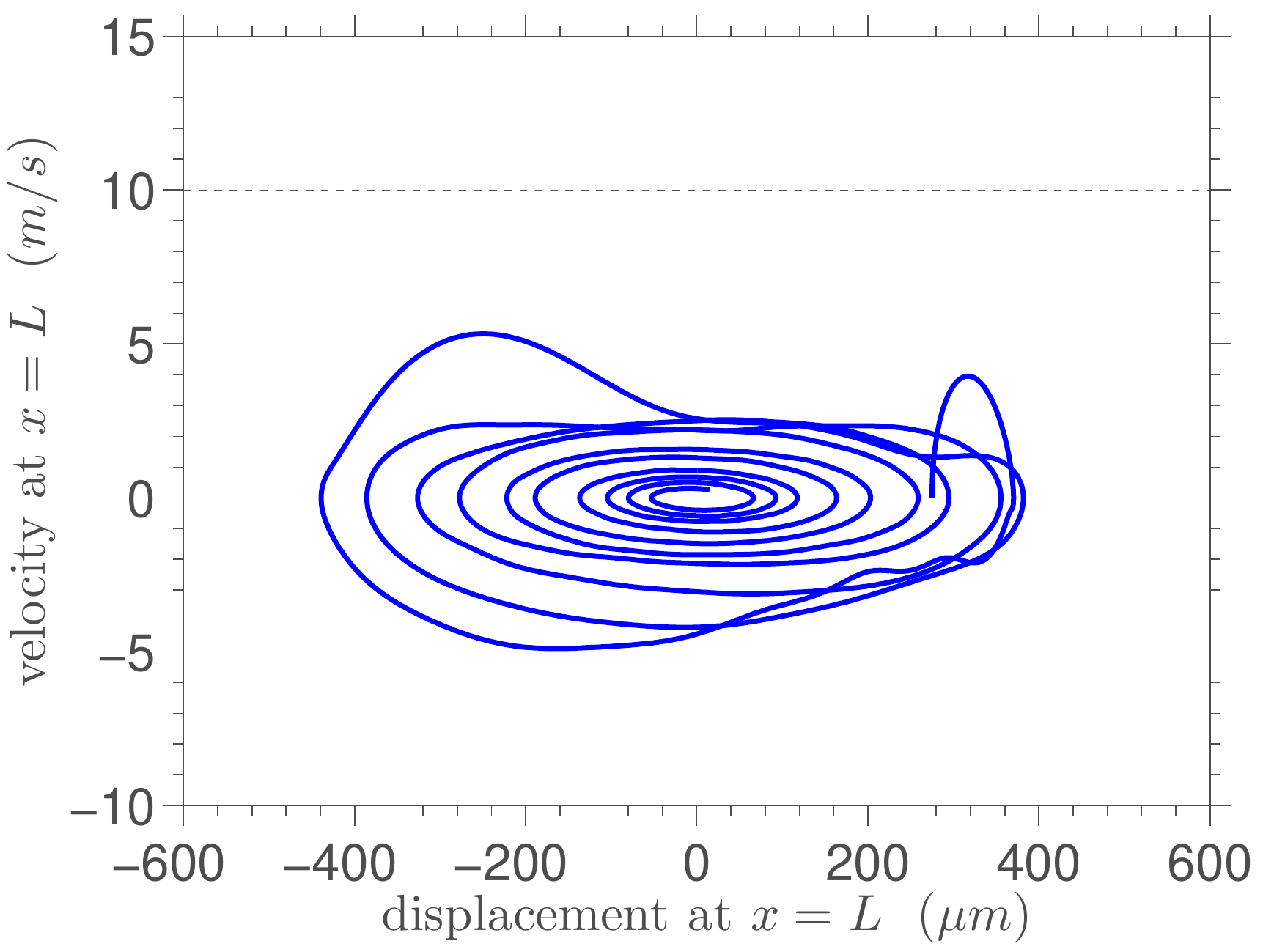}}\\
				\subfigure[$\rho A L/m \approx 1$]{
				\includegraphics[scale=0.38]{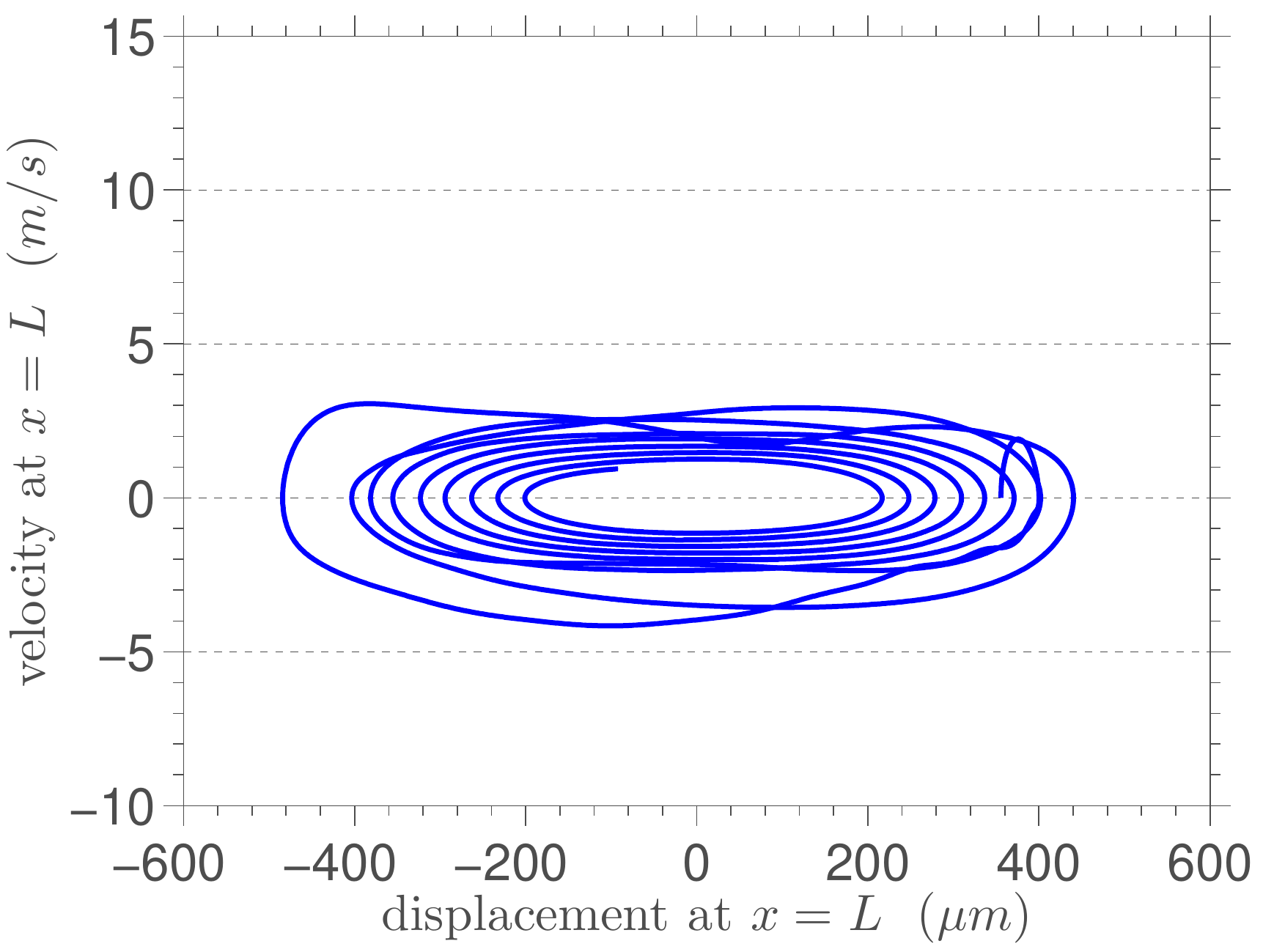}}
				\subfigure[$\rho A L/m \approx 0.2$]{
				\includegraphics[scale=0.38]{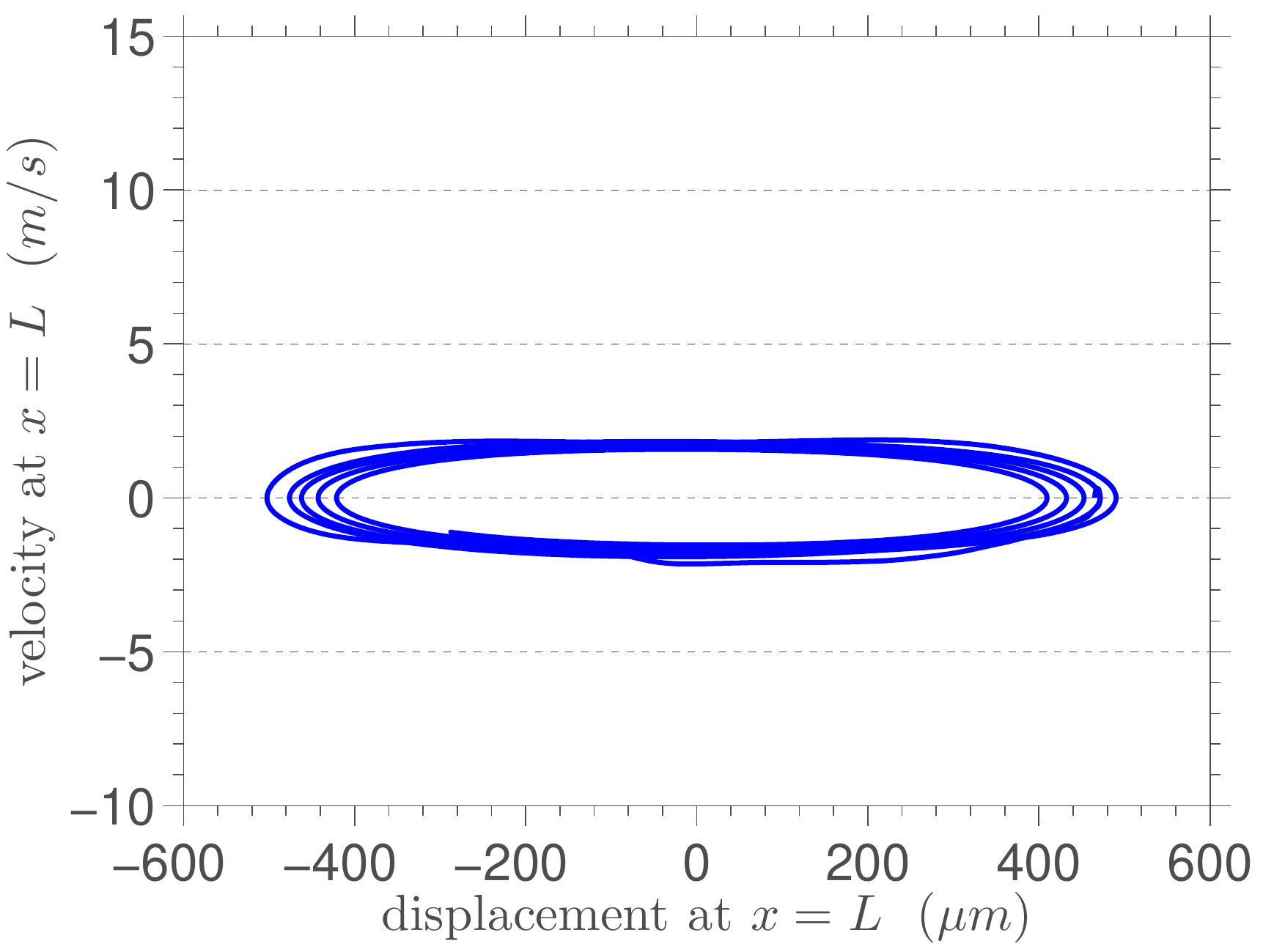}}
				\caption{This figure illustrates the mean value of the fixed-mass-spring bar
							 phase space at $x=L$, for several values of the continuous--discrete mass ratio.}
				\label{phase_space_fig}
\end{figure}

The mean phase space of the fixed-mass-spring bar at $x=L$ is shown, for different values 
of the continuous--discrete mass ratio, in Figure~\ref{phase_space_fig}. The observation made
in the previous paragraph can be confirmed by analyzing this figure, since the system mean
orbit tends from a stable focus to an ellipse as the continuous--discrete mass ratio decreases.
In other words the limiting behavior of the system when $\rho A L/m \to 0^+$ is a mass-spring
system.

\subsection{Analysis of Random System PDF}
\label{analysos_pdf}

\begin{figure} [ht!]
				\centering
				\subfigure[$\rho A L/m \approx 10$]{
				\includegraphics[scale=0.38]{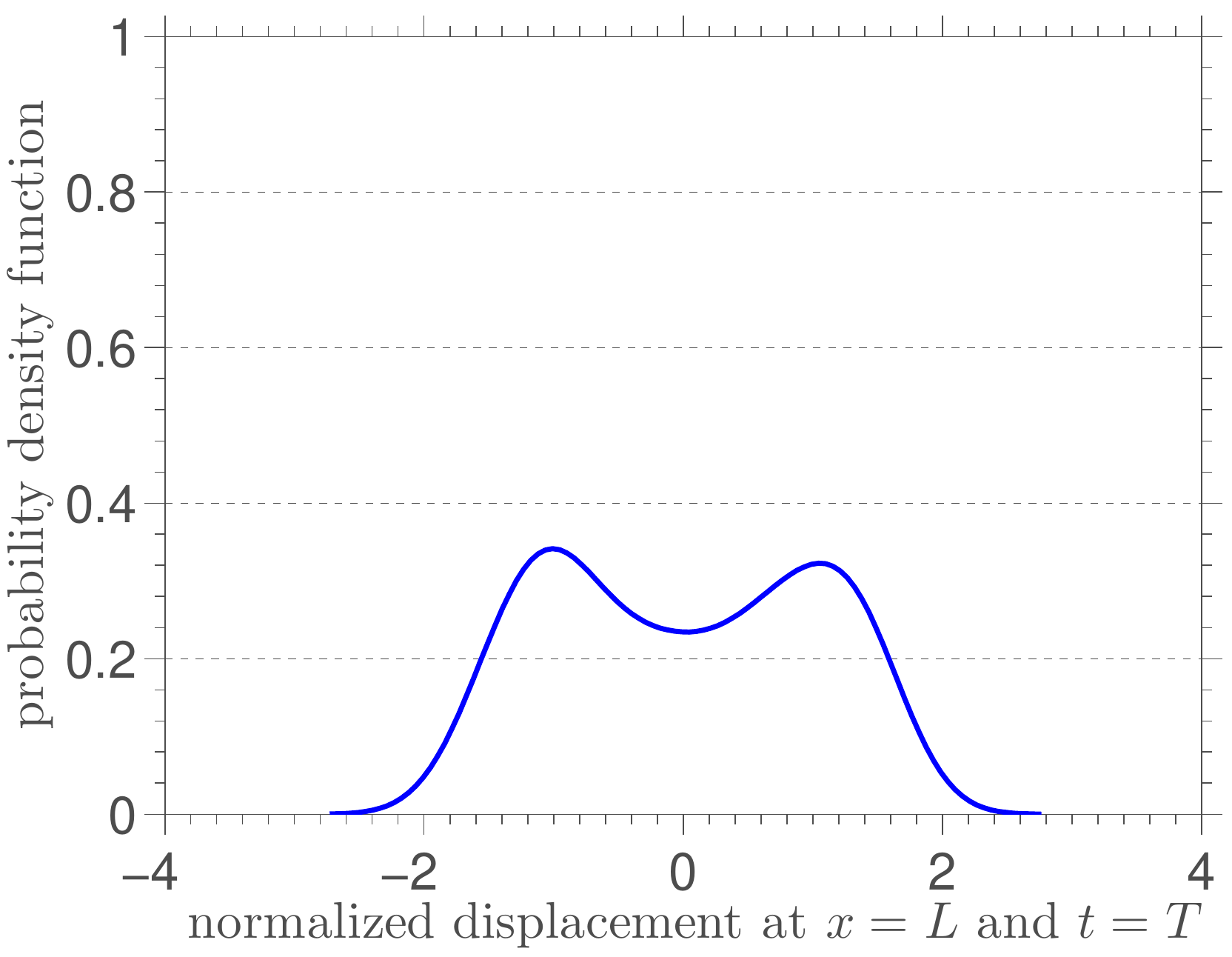}}
				\subfigure[$\rho A L/m \approx 2$]{
				\includegraphics[scale=0.38]{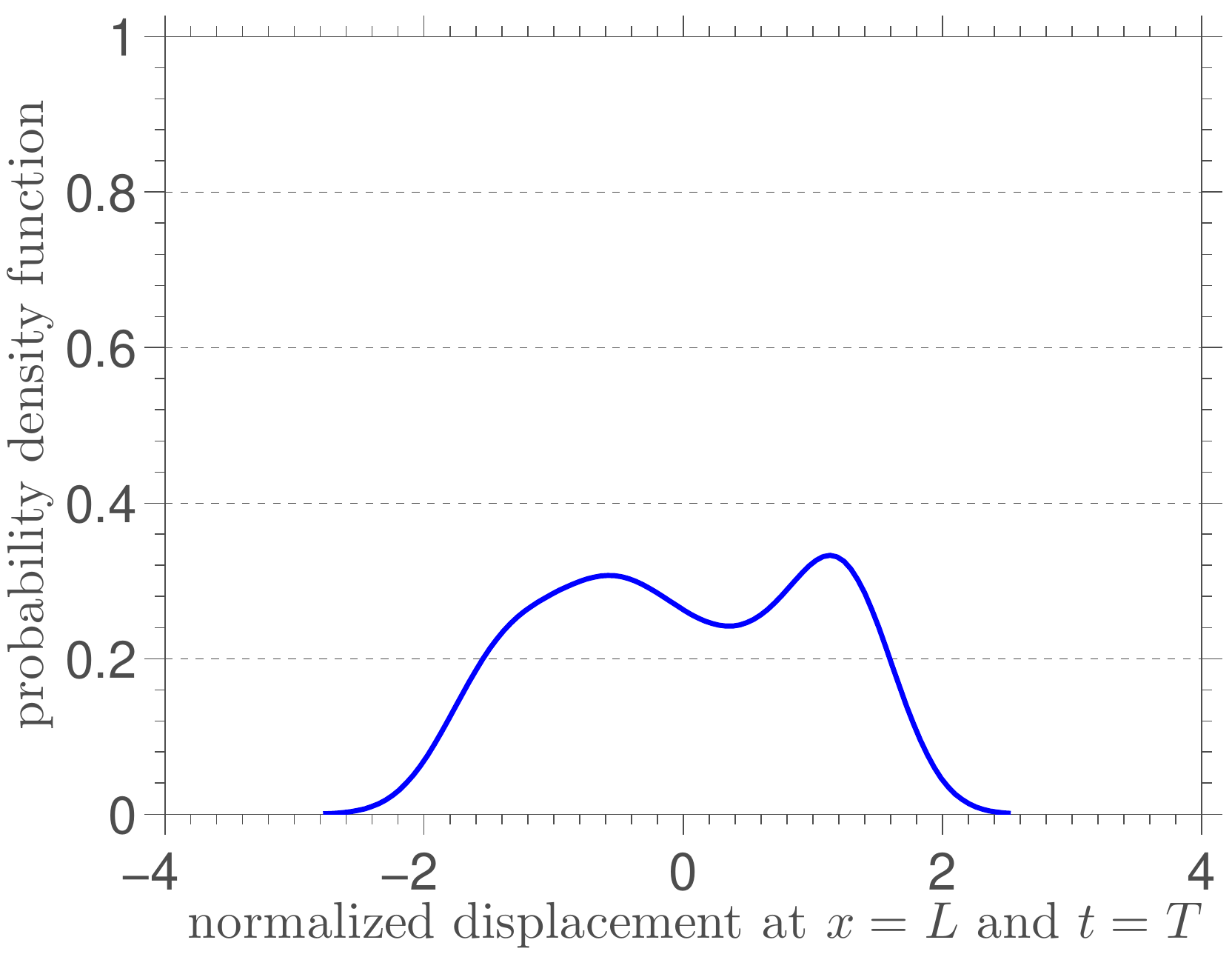}}\\
				\subfigure[$\rho A L/m \approx 1$]{
				\includegraphics[scale=0.38]{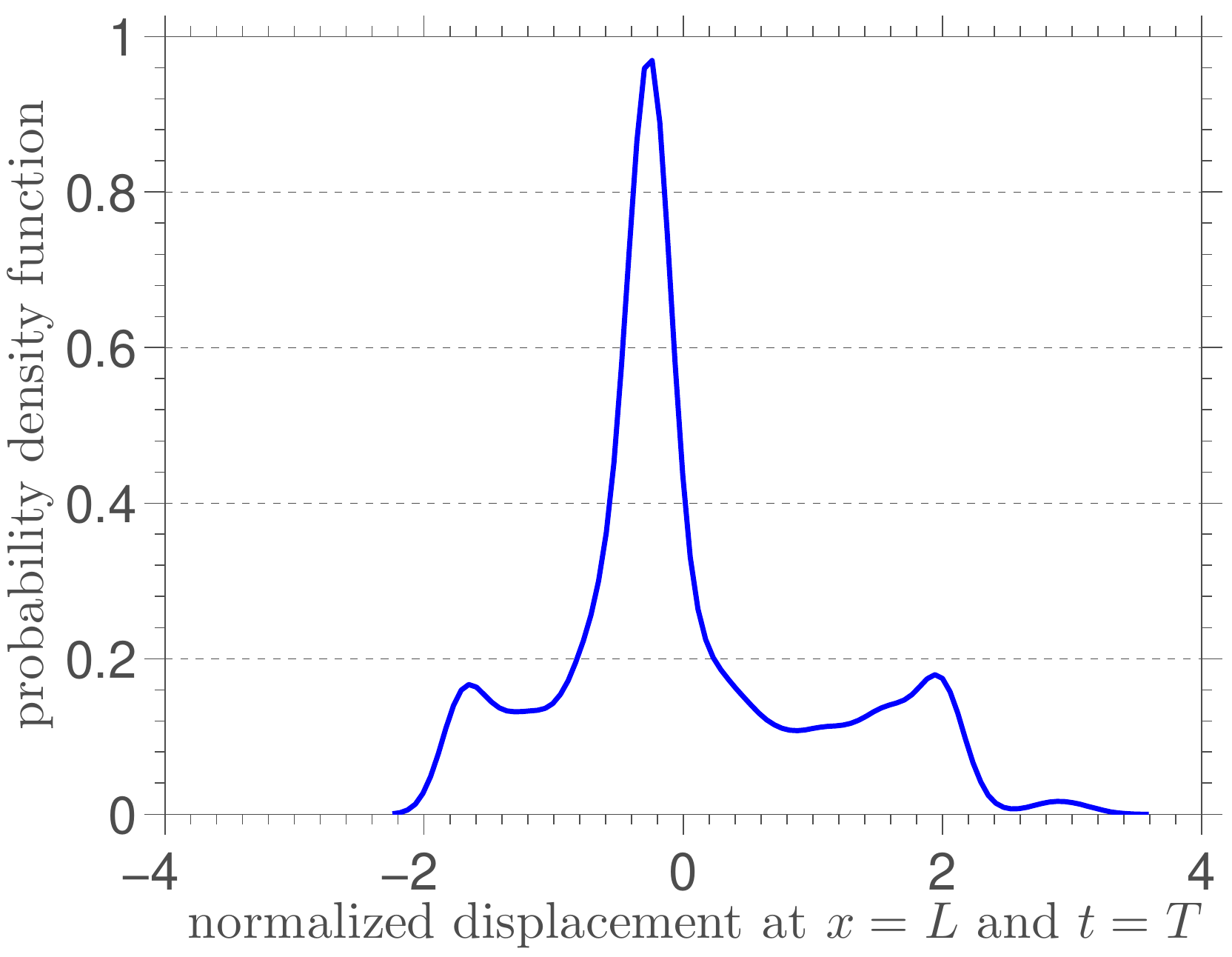}}
				\subfigure[$\rho A L/m \approx 0.2$]{
				\includegraphics[scale=0.38]{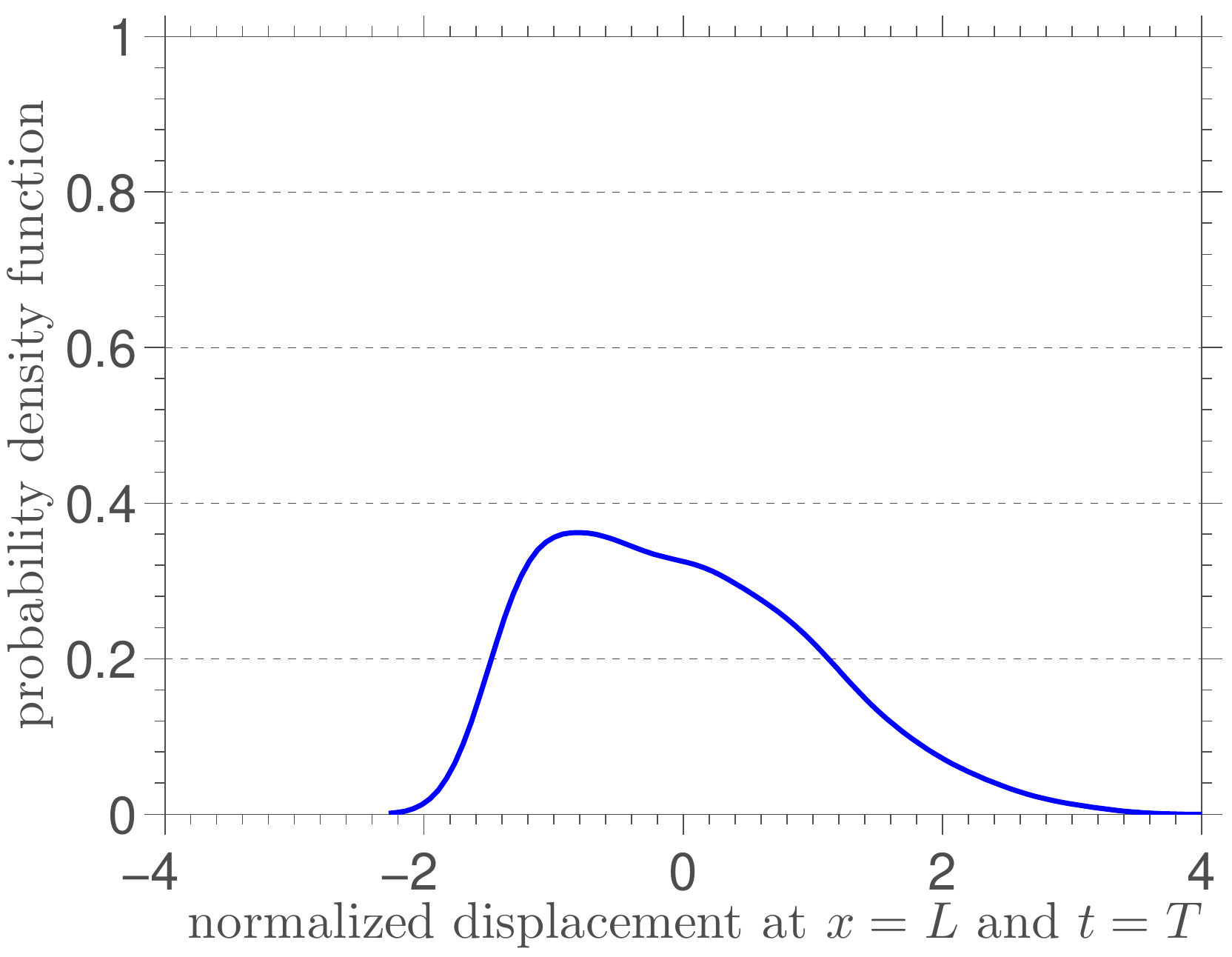}}
				\caption{This figure illustrates estimations to the PDF of the (normalized)
								random variable $\randproc{U}(L,T,\cdot)$, for several values 
								of the continuous--discrete mass ratio.}
				\label{pdf_uL_fig}
\end{figure}

The difference between the system dynamical behavior, for different values
of $m$, is even clearer if one looks to the PDF estimations of the (normalized)
random variable $\randproc{U}(L,T,\cdot)$, which are presented in Figure~\ref{pdf_uL_fig}.
For large values of the continuous--discrete mass ratio, the PDF
of $\randproc{U}(L,T,\cdot)$ displays bimodal shape, which tends to
a unimodal shape as the lumped mass grows, i.e., the continuous--discrete 
mass ratio decreases.
Furthermore, it can be noted that when $\rho A L/m \approx 1$
the greatest probability occurs around the mean value of 
$\randproc{U}(L,T,\cdot)$.

\section{CONCLUDING REMARKS}
\label{concl_remaks}

A model to describe the dynamics of fixed-mass-spring bar
with a random elastic modulus is presented. The aleatory parameter 
is modeled as a random variable with gamma distribution, 
being the probability distribution of this parameter obtained by 
the principle of maximum entropy. The paper analyzes some 
configurations of the model to order to characterize the effect 
of the lumped mass in the overall behavior of this dynamical system.
This analysis shows that the dynamics of the random system is 
significantly altered when the values of the lumped mass are varied.


\section{ACKNOWLEDGMENTS}

The authors are indebted to Brazilian Council for Scientific and 
Technological Development (CNPq), Coordination of Improvement 
of Higher Education Personnel (CAPES), and Foundation for Research 
Support in Rio de Janeiro State (FAPERJ) for the financial support.

\newpage


\end{document}